# Reflection Principles for C∞ Harmonic Maps between Riemannian Manifolds and a Schwarz Reflection Principle for Harmonic and Holomorphic Maps from a Class of Hermitian Symmetric Spaces

Dominic S.P. Leung

To Samuel 俊智 and Elliott 俊贤

**Abstract.** After establishing the uniqueness of the continuation of local Cauchy data for C∞ harmonic maps between two Riemannian manifolds M and N, we prove (i) a reflection principle for a C∞ minimal submanifold Y of a Riemannian manifold M that contains a reflective submanifold of M as a hypersurface and (ii) the reflection principle of for a class of C∞ harmonic maps between Riemannian manifolds that admit reflective submanifolds. A reflection principle for harmonic functions on M with respect to a reflective hypersurface of M follows easily from the above reflection principle for harmonic maps. To consider reflection principles for harmonic and holomorphic maps from general Hermitian symmetric spaces M, the notion of recursive real forms is introduced. Hermitian symmetric spaces admitting such recursive real forms include Hermitian manifolds of constant holomorphic curvatures, the complex hyper-quadrics and their respective non-compact duals. Let M be an irreducible Hermitian symmetric space of complex dimension n with a recursive real form $B_1$ and corresponding involutive isometry $\sigma_1$, X be a C∞ Riemannian manifold with $B_2$ as a reflective submanifold and corresponding involutive isometry $\sigma_2$. Also let h : M → X be a harmonic map such that (i) $h(B_1) \subset B_2$ and (ii) $h_*(w) \in T_{h(q)} B_2^\perp$ for $w \in T_q B_1^\perp$ for all q ϵ $B_1$. We will prove that $h(\sigma_1(p)) = \sigma_2(h(p))$ for all p ϵ M. A special case of the previous assertion is that X is a Hermitian symmetric space, $B_2$ is a real form of X and h is a holomorphic map. These two assertions generalize the classical Schwarz reflection principle for harmonic functions or holomorphic functions on the complex plane to harmonic maps and holomorphic maps from a class of Hermitian symmetric spaces admitting recursive real forms. Furthermore, a reflection principle for meromorphic functions defined on the complex plane and taking on real values on the real line is also proved.

*Key* **words and phrases**: Harmonic maps, linear sigma model, nonlinear sigma model, unique continuation of harmonic maps, Cauchy problem, reflective submanifolds, reflection principle, real forms, recursive real forms, holomorphic map, holomorphic functions, meromorphic functions, harmonic functions, Riemannian symmetric spaces, Hermitian symmetric spaces, involutive isometry, Schwarz reflection principle, bounded symmetric domains

**1. Introduction**

For simplicity, we will assume that all differential geometric structures to be C∞ unless otherwise stated. One of the purposes of this paper is to prove a reflection principle similar to that in [L2] for C∞ harmonic maps between Riemannian manifolds. Let M and N be two simply connected Riemannian manifolds. This paper begins with the establishing the uniqueness of the continuation of local Cauchy data for C∞ harmonic maps between M and N, making use of the classical results on the unique continuation of harmonic maps [AR2, SA, XY] and some elementary facts in PDE. Using these results and the notion of reflective submanifolds, we prove (i) a reflection principle for a $C^\infty$ minimal submanifold Y of a Riemannian manifold M that contains a reflective submanifold of M as a hypersurface and (ii) the reflection principle of for a class of $C^\infty$ harmonic maps between Riemannian manifolds that admit reflective submanifolds. A reflection principle for harmonic functions on symmetric space M with respect to a reflective hypersurface of M follows easily from the above reflection principle for harmonic maps. Recall that for a Hermitian symmetric space M the fixed-point sets of antiholomorphic involutive isometries (called complex conjugations) are called the real forms of M [L5]. To consider reflection principles for harmonic and holomorphic maps of general Hermitian symmetric spaces M, the notion of a recursive real form B of M is introduced. Recursive real forms include the real forms of Kahler manifolds of constant holomorphic curvatures and real forms of the complex hyperquadrics and their respective non-compact duals as irreducible Hermitian symmetric space. Let M be an irreducible Hermitian symmetric space of complex dimension n with a recursive real form $B_1$ and corresponding involutive isometry $\sigma_1$, X be a $C^\infty$ Riemannian manifold with $B_2$ as a reflective submanifold and corresponding involutive isometry $\sigma_2$. Also let h : M → X be a harmonic map such that (i) $h(B_1) \subset B_2$ and (ii) $h_*(w) \in T_{h(q)} B_2^\perp$ for w ∈ $T_q B_1^\perp$ for all q ∈ $B_1$. We will prove that $h(\sigma_1(p)) = \sigma_2(h(p))$ for all p ∈ M. A special case of the previous assertion is that X is a Hermitian symmetric space, $B_2$ is a real form of X and h is a holomorphic map. These two assertions generalize the classical Schwarz reflection principle for harmonic functions or holomorphic functions on the complex plane to harmonic maps and holomorphic maps from a class of Hermitian symmetric spaces admitting recursive real forms. Furthermore, a reflection principle for meromorphic functions defined on the complex plane and taking on real values on the real line is also established.

Part of the research for this paper was performed while the author was a visiting scholar at the Mathematical Science Center at Tsinghua University in the summer and fall of 2012. The author would like to express his thanks for the financial and facility supports, hospitality and the encouragements from S.T. Yau, Stephen S.T. Yau and Y.S. Poon received during the visit. In particular, he would like to thank S.T. Yau for his invitation for the visit. The author would also express his thanks to P.T. Ho of University of Maryland for providing copies of some very old references related to PDE used in this paper.

**2. Unique continuation of local Cauchy data for harmonic maps between Riemannian manifolds.**

The study of harmonic maps began with pioneer work of Ells-Simpson [ES]. Let *M* and *N* be smooth oriented Riemannian manifolds of dimension m and n respectively with metrics $ds_M^2$ and $ds_N^2$. Let *h: M → N* be a smooth map. Then the inverse image, $h^*ds_N^2$, considered as a form on *M*, is positive semi-definite and $Tr(h^*ds_N^2)$, its trace relative to $ds_M^2$, is a function on M taking non-negative values. If D is a compact domain in M, the integral

(2.1) $$E(h,D) = \frac{1}{2} \int_D Tr(h^*ds_N^2) * 1$$

where $*1$ is the volume element of *M*, is called the energy of *h* over *D*. The map *h* is called harmonic of over *D* if the energy *E(h,D)* has a critical value relative to all map which agree with *h* on the boundary of *D*. *h* is called a harmonic map if it is harmonic over any compact domain *D* on *M*. Let X and Y be cross sections of *TM* and *dh* be the differential 1-form with values in the pulled back bundle $f^{-1}TN$ defined by

(2.2) $$dh\,(X) \stackrel{\text{def}}{=} h_*(X)$$

*dh* is a cross section of the vector bundle $TM^* \otimes f^{-1}TN$. Let $\nabla$ denotes the induced connection on the bundle $TM^* \otimes f^{-1}TN$ from the natural Riemannian connections on *M* and *N*. Then the tension field of the harmonic map h can be defined by for any orthonormal local frame fields $e_i$

(2.3) $$\tau(h) \stackrel{\text{def}}{=} (\nabla_{e_i} dh) e_i$$

It can be proved in [J1] and [XY] that h is an harmonic map if and only if $\tau(h) \equiv 0$. Let U be a local coordinate neighborhood with coordinate systems $\{x^i\}$ and $\{y^\alpha\}$ of *M* and *N* respectively the condition $\tau(h) \equiv 0$ can be written as ( [J1] and [XY])

(2.4) $$\Delta_M h^\alpha + g^{ij} \bar{\Gamma}^\alpha_{\beta\gamma} \frac{\partial h^\beta}{\partial x^i} \frac{\partial h^\gamma}{\partial x^j} = 0,$$

where $\Delta_M$ and $g^{ij}$ denote respectively the Laplacian operator and metric tensor of M, $\bar{\Gamma}^\alpha_{\beta\gamma}$ the Christofell symbols of N, $1 \leq i, j \leq m$ and $1 \leq \alpha, \beta, \gamma \leq n$. Let *B* be a piece of hypersurface in an open set *U*, a set of Cauchy data for the system of partial differential equations (2.4) is given by the Jacobian matrix of partial derivatives

(2.5) $$\left[\frac{\partial h^\beta}{\partial x^i}\right].$$

The matrix (2.5) also defines for every p ϵ B, the tangential map, $h_*: T_pM \to T_{h(p)}N$. Conversely tangential map $h_*$ of a harmonic map defines, when restricted to *B* , a set of Cauchy data on *B* for the system of partial differential equations (2.4).The facts in the following theorem are proved in [AN2]

**Theorem 2.1.** *Let A be a linear elliptic second-order differential operator on a domain D of $R^n$. In D let u = ($u^1$, ..., $u^n$) be functions satisfying the differential inequalities*

$$|Au| \leq constant \left\{\sum_{i,\beta} \left|\frac{\partial u^\beta}{\partial x^i}\right| + \sum_\beta |u^\beta|\right\}.$$

*Then if u either vanishes to an infinite order at a single point in D or u is zero in an open subset, then u is identically zero in D.*

Using the results in theorem above the following facts are proved in [SA]

**Theorem 2.2**. *Let $f^1$, $f^2$ be two harmonic maps $M \to Y$ and M connected. If either they agree to an infinite order at some point or they agree on an open subset of M, then they are identical*

With these preparations, we will prove the following main results of this section of the paper.

**Theorem 2.3.** *Let $f^1$, $f^2$ be two harmonic maps $M \to Y$ and M connected. If they satisfy the same Cauchy data, or equivalently their differentials have the same values, along a piece of hypersurface B of M, then they are identical. In particular, if, $f = f^1 = f^2$, and $f\mid_B$ is an imbedding, then the submanifold f(M) of Y is uniquely determined by the distribution of tangent spaces, $f_*(T_q M)$ for $q \in B$, along B.*

Proof: Since $f^1(B) = f^2(B)$, we can assume without loss of generality that B lies in a coordinate neighborhood of M with local coordinates $\{x^i\}$ and their images lie inside a coordinate neighbors of Y with local coordinates $\{y^i\}$. Then the local expressions, $f^{1\beta}$ and $f^{2\beta}$ of $f^1$ and $f^2$ (resp.) satisfy a system of quasi-linear 2$^{nd}$ order elliptic PDE of the form (2.4). As such the hypersurface *B* is a non-characteristic Cauchy hypersurface of (2.4). It is an elementary fact in PDE that, at any point of *B*, all the partial derivatives of $f^{1\beta}$ and $f^{2\beta}$ with respect to $x^i$ are completely determined by their Cauchy data and coefficients of the system of PDE (2.4) (see for example [EL] or [JF]). Therefore $f^{1\beta}$ and $f^{2\beta}$ agree to an infinite order at all points of B. It follows from Theorem 2.2 that $f^1$ and $f^2$ are identical □

**Remark 2.5** Theorems 2.1 and 2.2, and consequently Theorem 2.3 are true under much weaker assumptions than $C^\infty$ ([AN2] and [KJ]). But the $C^\infty$ assumption is adequate for our current applications. The uniqueness result of the solutions of the system (2.4) for a given Cauchy data set is well known, in fact it was pointed out in [AN2] that such uniqueness results could follow easily from the main results in that paper [AN2] for a general class of quasi-linear 2$^{nd}$ order elliptic systems.

**Remark 2.6** The results of [AN2] were announced in [AN1], although the reference [AN2] is no longer available in U.S. libraries according to the search engine WorldCat, but its preprint [AN3] is available in many libraries in U.S.

**Remark 2.7** The existence of solutions for given elliptic system of partial differential equations satisfying a given set of Cauchy data is not well-posed even for linear systems in the $C^\infty$ category. However, under the assumption of analyticity, the existence and uniqueness of harmonic maps and also of minimal submanifolds can be proved using the classical Cartan Kähler theory ([K] and [BCG]) by developing appropriate exterior differential systems in appropriate bundle spaces; these have been carried out in [L6] and [L1] respectively. The approach used herein gives a more streamline and simpler proof that is more accessible to geometers and analyst without any exposure to Cartan Kähler theory.

In the remainder of this paper, we will apply Theorem 2.3 to study reflection principles in differential geometry.

## 3. Application to a reflection principle for minimal submanifolds

Under the assumption of real analyticity, a reflection principle for minimal submanifolds has been proved in [L2]. We will give herein a proof the same reflection principle for $C^\infty$ minimal submanifolds. Recall that a reflective submanifold [L2] B of a Riemannian manifold N, is a connected component of an involutive isometry $\rho_B$ of N.

**Theorem 3.1**. *Let h : M → N be a minimal submanifold of M, U and V be coordinate neighborhoods of M, N respective such that, $h|_U$ is an imbedding, $h(U) \subset V$ and B be a reflective submanifolds of N with defining involutive isometry $\rho_B$. If $h(U)$ contains $B \cap V$ as a hypersurface, then $\rho_B(h(M)) = h(M)$.*

Proof: It is well known [CG] that h defines a minimal immersion of M into N, if and only if h is an isometric harmonic map. We can assume without loss of generality the both U and V are local coordinate neighborhoods of M and N respectively. Therefore, in local coordinate systems $\{x^i\}$ and $\{y^\alpha\}$ of M and N respectively, h satisfies the second order quasi-linear systems (2.4). For any point $p \in U$ and $w \in T_{h(p)} B^\perp$, then $T_{h(p)} h(M) = T_{h(p)} B \oplus W$, with W being the subspace of $T_{h(p)} N$ spanned by w. Since $\rho_{B*}(w) = -w$, $\rho_B$ leaves the $T_{h(p)} h(M)$ invariant. Therefor the two harmonic maps h and $\rho_B \circ h$ have the same Cauchy data along $B \cap V$. By Theorem 2.2 we have $h = \rho_B \circ h$. □

**Remark 3.1** Recall that according to the classification of reflective submanifolds of Riemannian symmetric spaces ([L3] and lL4]), each irreducible symmetric space contains many non-isometric reflective submanifolds. It has been proved in [G] that if B is a reflective submanifolds of a irreducible noncompact symmetric space N, then for any point $p \in B$ and $w \in T_p B^\perp$, B can be extended to a minimal submanifold, $f: \mathbf{R} \times B \to N$, of N such that $f(0, p) = p$ for all $p \in B$ and $\frac{\partial f(t,p)}{\partial t}\big|_{t=0} = w$. In fact, this is accomplished by action on *B* of a properly chosen one parameter family of isometries of *N*. The methodology established in [G] to enlarge totally geodesic submanifiolds of *N* to minimal submanifolds of one dimension higher can be used for a large class of geodesic submanifolds that contains the reflective submanifolds as a proper subset. Hence there is a very large number of examples of known minimal submanifolds for which Theorem 3.1 can be applied. Theorem 3.1 does not require that the minimal submanifold *f(M)* to be real analytic, there it should be of interest to find examples of non-analytic submanifolds for which the theorem could be applied.

The class of reflective submanifolds introduced in [L2] for Riemannian symmetric spaces is currently an active area of research see for examples [BENT], [G], [MQ] and [T].

**4 A reflection principle for harmonic maps.**

**Theorem 4.1.** *Let h : M → N be a harmonic map, $B_1$ and $B_2$ reflective submanifolds of M and N (resp.) with involutive isometries $\sigma_1$ and $\sigma_2$ (resp.) such that $h(B_1) \subset B_2$. Let H be a codimension one submanifold of M, H contains $B_1$ and is left invariant by $\sigma_1$. Furthermore,* if we also have,

(1) $h(\sigma_1(p)) = \sigma_2(h(p))$ for all $p \in H$,
(2) $h_* \circ \sigma_1(v) = \sigma_{2*} \circ h_*(v)$ for all $v \in T_p H$ and $p \in H$
(3) $h_*(w) \in T_{h(q)} B_2^\perp$ for $w \in T_q B_1^\perp$ for all $q \in B_1$.

*Then we have $h(\sigma_1(p)) = \sigma_2(h(p))$ for all $p \in M$.*

Proof: Consider the map $\tilde{h}: M \to N$ define by $\tilde{h} = \sigma_2 \circ h \circ \sigma_1$. Let $p \in H$, choose a local orthonormal frame fields $\{e_i\}$ in a neighborhood $\mathcal{U}$ of p, such that $\{e_1, ..., e_{m-1}\}$ spans $T_q H$ for all $q \in \mathcal{U} \cap H$. It follows from the assumption (3) that for $q \in B_1$ $h_*(e_m(q)) \in T_{f(q)} B_2^\perp$. Therefore $\sigma_{2*}(h_*(e_m(q))) = - h_*(e_m(q))$. It follows that $\tilde{h}_*(e_m(q)) = h_*(e_m(q))$. Therefore, it follows from assumption (2) that $\tilde{h}_*(v) = h_*(v)$ for all $v \in T_q M$ with $q \in \mathcal{U} \cap B_1$. $\tilde{h}$, being also a harmonic map, it follows from Theorem 2.3 that $\tilde{h}$ and $h$ are identical on M. □

**Remark 4.1.**

(i)   If in Theorem 4.1, $H = B_1$, then the conditions (1) and (2) are automatically satisfied.
(ii)  For a general harmonic map condition (3) is in general not true, additional assumptions such as isometry, conformality, holomorphicity, or others will be needed to force this to be true.
(iii) Harmonic maps are known as the nonlinear sigma model in physics [J2], it is not clear if there will be any applications of the conclusions of Theorem 2.3 and Theorem 4.1 to problems in physics.

**Remark 4.2.** All the examples of reflective submanifolds known to the author are analytic submanifolds of symmetric spaces ([L4] and [L5]). In view of the fact that Theorem 4.1 holds true for reflective submanifolds $B_1$ and $B_2$ that need not be analytic, there it will be of interest to find examples of reflective submanifolds that are not analytic.

**5. Application to a reflection principle for harmonic functions with respect to reflective hypersurfaces**
When N is the real line, the energy integral is nothing but the classical Dirichlet integral, the torsion field is just the Laplacian of a function on M. This is also called the linear sigma model in physics [J2].

In the special cases, when M is a space of constant curvatures, N is the real line, $B_1$ (= H) a codimension one totally geodesic submanifold and $B_2$ the 0 point, Theorem 4.1 specializes to the following theorem.
**Theorem 5.1**. *Let g be a harmonic function on a simply connected space K, of constant curvature $k \leq 0$, that takes on the value 0 on codimension one totally geodesic submanifold B with geodesic reflection $\rho$, then $g(\rho(x)) = -g(x)$ for all $x \in M$.*

For harmonic functions on the flat Euclidean spaces $E^n$, Theorem 5.1 is well known, see for example [ABR] for a proof using a completely different approach.

**6. Generalization of the Schwarz reflection principle to harmonic maps and holomorphic maps defined on a class of higher dimension Hermitian symmetric spaces and meromorphic functions defined on the complex line.**

Recall that the real form of simply connected Hermitian symmetric space M is the fixed-point set of an involutive anti-holomorphic isometry of M. Such real forms for irreducible Hermitian symmetric spaces have been completely classified [L5], in particular the classical bounded symmetric domains. For completeness and easy references, a table of all such real forms of classical bounded symmetric domains is included in the appendix of this paper. It is well known that holomorphic maps are examples of harmonic maps ([ES] and [ XY])

**Lemma 6.1.** *Under the assumptions of Theorem 4.1, if furthermore M and N are Hermitian symmetric spaces, h is a Holomorphic map and $B_1$ and $B_2$ real forms M and N respectively, then condition (3) of Theorem 4.1 is always satisfied.*

Prrof: Denote the complex structures of M and N by $J_M$ and $J_N$ respectively. For *all $q \in B_1$ and any w $\in T_q B_1^{\perp}$*, there exists a $v \in T_q B_1$ such that $w = J_{M*}(v)$. Therefore, we have

$$h_*(w) = h_*(J_{M*}(v))$$
$$= J_{N*}(h_*(v)).$$

Since $h_*(v)\ T_{h(q)} B_2$, $J_{N*}(h_*(v)\ T_{h(q)} B_2^{\perp}$. □

Since every holomorphic map between Hermitian manifolds is a harmonic map [XY], in view of Lemma 6.1 specializing Theorem 4.1 to holomorphic maps we have the following theorem.

**Theorem 6.1.** *Let $f : M \to N$ be a holomorphic map between the two Hermitian symmetric spaces, $B_1$ and $B_2$ be the (resp.) real forms of M and N with involutive isometries $\sigma_1$ and $\sigma_2$ (resp.) such that $h(B_1) \subset B_2$. Let H be a real hypersurface of M, H contains $B_1$ and is left invariant by $\sigma_1$. Furthermore we also have,*

(1) *$f(\sigma_1(p)) = \sigma_2 (f(p))$ for all $p \in H$,*
(2) *$f_* \circ \sigma_1( v) = \sigma_{2*} \circ f_* (v)$ for all $v \in T_p H$ and $p \in H$*

*Then, we have $f(\sigma_1(p)) = \sigma_2 (f(p))$ for all $p \in M$.*

The above theorem naturally leads us to the following notion of recursive real forms.

**Definition 6.1**. *Let $\underline{M}^n$ be either $C^n$ or an irreducible Hermitian symmetric space of complex dimension n. A real forms B of $\underline{M}^n$ is called recursive if it satisfies the following conditions:*

*There exists a sequence of Hermitian symmetric spaces $\widetilde{M}_k$ of the same type,*

(6.1) $$\widetilde{M}_1 \subset \ ...\subset \widetilde{M}_k \subset\ ... \subset \widetilde{M}_{n-1} \subset M^n\ .$$

Such that

   (i)  $dim_c(\widetilde{M}_k) = k$,
   (ii) $\widetilde{M}_{n-1}$ ($\widetilde{M}_l$ ) *is a totally geodesic submanifold of $M^n$ (resp. $\widetilde{M}_{l+1}$),*
   (iii) *all the differential geometric structures of $\widetilde{M}_l$ are obtained by the restriction of the corresponding ones of $M^n$*
   (iv) *$B \cap \widetilde{M}_k$ is a real form of $\widetilde{M}_k$ for k =1, ..., l and $H_l = (B \cap \widetilde{M}_l) \cup \widetilde{M}_{l-1}$ is a smooth real hypersurface of $\widetilde{M}_l$ for $l = 2, ..., n$-1.*

**Remark 6.1.** *It follows from the above definition of recursive real forms that*
*(i) if f is holomorphic map of $M^n$ into another complex manifold N, then the restriction of f to $\widetilde{M}_l$ also defines a holomorphic map of $\widetilde{M}_l$ into N. This follow readily from the fact that the complex structure on $\widetilde{M}_l$ is simply the restriction of that of $M^n$ to its submanifold $\widetilde{M}_l$ and*
*(ii) if g is a real harmonic function on $M^n$, then the restriction of g to $\widetilde{M}_l$ also defines a real harmonic function on $\widetilde{M}_l$.*

**Proposition 6.1.** The following are examples of recursive real forms of Hermitian symmetric spaces:

|   | Type | Hermitian Symmetric Spaces | Real Forms |
|---|------|---------------------------|------------|
| 1 | Euclidean | $C^n$ | $R^n$ |
| 2 | Hermitian Hyperbolic | $SU(1, n)/S(U(1) \times U(n))$ | $SO(1, n)/SO(n)$ |
| 3 | Complex Projective Space | $SU(1+n)/S(U(1) \times U(n))$ | $SO(1+n)/SO(n)$ |
| 4 | Noncompact dual of Complex Hyperquadric | $SO(2, n)/SO(2) \times SO(n)$ | $[SO(1, q)/SO(1) \times SO(q)] \times [SO(1, n-q)/SO(1) \times SO(n-q)]$, $0 \le q < n/2$, n odd $0 \le q \le n/2$, n even > 2 |
| 5 | Complex Hyperquadric | $SO(2+n)/SO(2) \times SO(n)$ | $[SO(1+q)/SO(1) \times SO(q)] \times [SO(1+n-q)/SO(1) \times SO(n-q)]$, $0 \le q < n/2$, n odd $0 \le q \le n/2$, n even > 2 |

**Proof:**

The existence of sequence of submanifolds satisfying the conditions of Definition 4.1 for such a pair ($M^n$, B) of Euclidean type is obvious. The existence of such sequence for a pair ($M^n$, B) of complex hyperbolic type or BD I type is also not difficult to see, once we will realize the above two irreducible Hermitian symmetric spaces as bounded symmetric bounded domains for example as in [H]. We will denote by $M_{p,q}(R)$ and $M_{p,q}(C)$ the set of all p x q real and complex matrices respectively. We will examine each of the three cases above separately.

**Euclidean Type.** Let B be the real form $R^n$ naturally imbedded in $C^n$ and $\widetilde{M}_k = C^k$. In fact with respect to the global coordinates $z_i = x_i + \sqrt{-1}\, y_i$, $i = 1,\ldots, n$, the real hypersurface $H_l$ of $\widetilde{M}_l$ is defined by the equation $y_l = 0$ and B is defined by the equations $y_k = 0$ for $k = 1, \ldots, n$.

**Hermitian Hyperbolic Type.** In this case we can realize $M^n$ as the following bounded domain [H]:

(6.2) $\qquad \{Z \in M_{n,1}(C): {}^t Z \bar{Z} < 1\}$.

Let B be the unique real form of $M^n$, the n-dim real hyperbolic spaces, can be naturally identified with the n x 1 real matrices $M_{n,1}(R) \subset M_{n,1}(C)$ in a natural way, and $\widetilde{M}_k = \{Z \in M_{k,1}(C): {}^t Z \bar{Z} < 1\}$. Similarly with respect to the global coordinates $z_i = x_i + \sqrt{-1}\, y_i$, $i = 1, \ldots, n$, the complex real hypersurface $H_l$ of $\widetilde{M}_l$ is defined by the equation $y_l = 0$ and B is defined by the equations $y_k = 0$ for $k = 1, \ldots, n$.

**Complex Projective Space.** We will follow the notions for the complex projective space $P_n(C)$ as described in section 10.5 of [KN] closely. Let $z = [z^0: z^1: \ldots : z^n]$ be the homogeneous coordinates of $P_n(C)$ as the Grassmannian manifold of complex lines in $C^{n+1}$. In the coordinate neighborhood

$U_0 = \{ z; z^0 \neq 0\}$ of $P_n(C)$ consider the inhomogenous complex coordinates $(w^1, \ldots, w^n)$ with $w^k = z^k/z^0$ for $1 \le k \le n$. In this case we can realize $M^n$ as $P_n(C)$,

(6.3) $\qquad \widetilde{M}_k \;(=\; P_k(C)\;) = \{\, z \in P_n(C): z^{k+1} = 0, \ldots, z^n = 0\}$ for $k = n-1, \ldots, 1$,

and let B be the unique real form $P_n(R)$ of $P_n(C)$. The complex real hypersurface $H_k$ of $\widetilde{M}_k$ is defined by the equation $v^k = 0$ where we have set $w^k = u^k + \sqrt{-1}\, v^k$.

**Noncompact duals of Complex Hyperquadrics.** In this case we can realize the Hermitian symmetric space as the following bounded domain [H]:

(6.4) $\qquad \{X \in M_{2,n}(R): X^t X < I_2 \}$.

The complex structure $J_0$ in this case is given by

(6.5) $\qquad J_0 \begin{pmatrix} X_1 \\ X_2 \end{pmatrix} = \begin{pmatrix} X_2 \\ -X_1 \end{pmatrix}$.

For an n-dim row vector $X = (x_1, \ldots, x_n)$ and an integer q, $1 \le q \le n$, define the maps $p_{b,q}$ and $p_{e,n-q}$ from $R^n$ into $R^n$ by:

- $p_{b,q}(X) = (x_1, \ldots, x_q, 0, \ldots, 0)$ - retaining the 1st q components and setting the remaining ones to zero
- $p_{e,n-q}(X) = (0, \ldots, 0, x_{q+1}, \ldots, x_n)$ - retaining the last n - q components and setting the remaining ones to zero.

Note that we have $X = p_{b,q}(X) + p_{e,n-q}(X)$. Define the involutive map $\sigma_q$ of $M_{2,n}(R)$ by

(6.6) $\qquad \sigma_q \begin{pmatrix} X_1 \\ X_2 \end{pmatrix} = \begin{pmatrix} -p_{b,q}(X_1) + p_{e,n-q}(X_1) \\ p_{b,q}(X_2) - p_{e,n-q}(X_2) \end{pmatrix}$.

One can easily verify that $\sigma_q \circ J_0 = -J_0 \circ \sigma_q$. Hence $\sigma_q$ defines an anti-holomorphic isometry of the n-dim non-compact dual of the complex hyperquadric, as Hermitain symmetric space under consideration. It can easily be verified that $\sigma_q$ has the fixed point sets $B_q$, for $0 \le q \le [n/2]$, defined below:

(6.7) $\qquad B_q = \left\{ \begin{pmatrix} X_1 \\ X_2 \end{pmatrix} \Big| p_{b,q}(X_1) = 0 \text{ and } p_{e,n-q}(X_2) = 0 \right\}$

as the real forms of the noncompact Hermitian symmetric space under consideration. In this case we can set $\widetilde{M}_k = \{X \in M_{2,k}(R): X^t X < I_2\}$ and B be any one of the real forms $B_q$. With $X_1$ and $X_2$ as used in (6.5), we set $X_i = (x_{i1}, ..., x_{in})$, for $i = 1, 2$. With B set equal to $B_q$ then for $l = 1, ..., n$, the complex real hypersurface $H_l$ of $\widetilde{M}_l$ is defined by $x_{1l} = 0$, for $l \leq q$ and $x_{2l} = 0$ for $l > q$.

**Complex Hyperquadric.** We will follow closely the notions for the complex hyperquadric $Q_n(C)$ as described in Example 10.6 of [KN] wherein the isometry between $Q_n(C)$ and the real Grassmannian of 2-planes in $R^{n+2}$, $G_{2,n+2}(R)$, as a Hermitian symmetric space is also described. We will use the results in Example 10.6 of [KN] freely in the following discussions. In terms of the homogeneous coordinates $z = [z^0: z^1: ...: z^{n+1}]$ of $P_{n+1}(C)$, $Q_n(C)$ can be defined by the equations

(6.8a) $\qquad (z^0)^2 + (z^1)^2 + ... + (z^{n+1})^2 = 0.$

If we write the vector $z$ in $C^{n+2}$ as

(6.8a) $\qquad z = x + \sqrt{-1}\, y, \quad x = (x^0, x^1, ..., x^{n+1}) \ \&\ y = (y^0, y^1, ..., y^{n+1}) \in R^{n+2},$

then we can equivalently define $Q_n(C)$ by

(6.8b) $\qquad Q_n(C) = \{z \in P_{n+1}(C): \langle x, y \rangle = 0 \text{ and } \|x\| = \|y\|\}.$

The description of the complex hyperquadric as a Hermitian symmetric space and its Cartan decomposition is given in details in Example 10.6 of [KN]. Thus by duality, we transfer all the information related the recursive real forms of its dual, at least in the Lie algebra level to $G_{2,n+2}(R)$. We will also give a description of the manifolds structure related the related recursive real forms as subsets of $Q_n(C) \subset P_{n+1}(C)$. Let $z = [z^0: z^1: ...: z^{n+1}]$ be the homogeneous coordinates of $P_{n+1}(C)$. In this case we can realize $M^n$ as $Q_n(C)$,

(6.9) $\qquad \widetilde{M}_k = Q_k(C) = \{z \in Q_n(C): z^{k+1} = 0, ..., z^{n+1} = 0\}$ for $k = n-1, ..., 1$.

It is also proved in Example 10.6 of [KN] that the real tangent space of at the point $\beta_0 = (1/\sqrt{2}, \sqrt{-1}/\sqrt{2}, 0, ..., 0)$ of $Q_n(C)$ is naturally the image of the tangent space of the Hermitian symmetric space $Q_n(C)$ under the Cartan decomposition and under the diffeomorphism described there in. Using this complex coordinate system of $P_{n+1}(C)$, we can identify the real forms of $Q_n(C)$ that has been classified in [L5] as a subset of $P_{n+1}(C)$. In the neighborhood

(6.12) $\qquad U_0 = \{z \in P_{n+1}(C): z^0 \neq 0\},$

we can set $w^k = z^k/z^0$, then (6.8) becomes

(6.13) $\qquad 1 + (w^1)^2 + (w^2)^2 ... + (w^{n+1})^2 = 0.$

The restrictions of the coordinate functions $\{w^2, ..., w^{n+2}\}$ to $Q_n(C)$ define a complex coordinate system of the two open subsets $V_1$ and $V_1$ of $Q_n(C)$ defined respective by by

(6.14a) $\qquad w^1 = \sqrt{-1 - (w^2)^2 - \cdots - (w^{n+1})^2}$ and $1 + (w^2)^2 + \cdots + (w^{n+1})^2 \neq 0,$

(6.14b) $\qquad w^1 = -\sqrt{-1 - (w^2)^2 - \cdots - (w^{n+1})^2}$ and $1 + (w^2)^2 + \cdots + (w^{n+1})^2 \neq 0.$

Put

(6.15a) $\zeta^k = w^{k+1}$, $\zeta^k = \chi^k + \sqrt{-1}\upsilon^k$ for k =1, ..., n

(6.15b) $\zeta = (\zeta^1, ..., \zeta^n)$, $\chi = (\chi^1, ..., \chi^n)$ and $\upsilon = (\upsilon^1, ..., \upsilon^n)$.

Note that points on the open sets $V_1$ and $V_1$ are mirror images of one and other with identical $\zeta$ coordinate values up to their common boundary defined by $w^1 = 0$. For $1 \le q \le n$ define a map $\tau_q$ of $V_i$ to $V_i$ for both i =1,2 by

(6.15) $\tau_q(\zeta) = -p_{b,q}(\chi) + p_{e,n-q}(\chi) + \sqrt{-1}\left(p_{b,q}(\upsilon) - p_{e,n-q}(\upsilon)\right)$, for $\zeta$ in $V_1$ or $V_2$,

with the projection operators $p_{b,q}$ and $p_{e,n-q}$ defined as above. One can readily verify that

(6.16) $\tau_q(\sqrt{-1}\zeta) = -\sqrt{-1}\tau_q(\zeta)$ for $\zeta$ in $V_1$ or $V_2$.

Therefor $\tau_q$ defines an anti-holomorphic map of $V_i$ to $V_i$ for both i =1,2. For each value of q, $1 \le q \le n$ and $1 \le i \le 2$ let

(6.17) $\tilde{B}_{q,i} = \{\zeta \in V_i : p_{b,q}(\chi) = 0 \text{ and } p_{e,n-q}(\upsilon) = 0\}$

and let $\tilde{B}_q$ be the closure of the union of $\tilde{B}_{q,1}$ and $\tilde{B}_{q,2}$ in $Q_n(C)$ that is

(6.18) $\tilde{B}_q = \overline{\tilde{B}_{q,1} \cup \tilde{B}_{q,2}}$.

It follows that $\tilde{B}_q$ defines a real form of $Q_n(C)$. Note that for k= n-1,..., 1, $\tilde{B}_q \cap Q_k(C)$ is also a real form of $Q_k(C)$. For i =1,2 , and $l$ = n, n-1, ..., 2, let $\widetilde{H_{l,\iota}}$ be subset of $V_i$ be defined by

(6.19) $\widetilde{H_{l,\iota}} = (\tilde{B}_{q,i} \cap \tilde{M}_l) \cup (\tilde{M}_{l-1} \cap V_i)$,

then $\widetilde{H_{l,\iota}}$ as a hypersurface of $\tilde{M}_l$ can be defined by $\chi^l = 0$, for $l \le q$ and $\upsilon^l = 0$ for $l > q$. Let $\widetilde{H_l}$ be the closure of the union of $H_{l,1}$ and $H_{l,2}$ in $\tilde{M}_l$. Then $\widetilde{H_l}$ is a real hypersurface of $\tilde{M}_l$. We have completed the construction of an recursive real form of $Q_n(C)$.

In terms of the homogenous coordinates z of $P_{n+1}(C)$. $\tilde{B}_q$ can also defined implicitly by (6.8a) and

(6.20) $x^2 =0, ..., x^{q+1} =0$ and $y^{q+2} =0, ..., y^{n+2} = 0$,

with $x^i$ and $y^j$ as defined in (6.8b). □

**Theorem 6.2.** *Let $\underline{M}^n$ be a Hermitian symmetric space of complex dimension n with a recursive real form $B_1$ and corresponding involutive isometry $\sigma_1$, and X be a Riemannian manifold with a reflective submanifold $B_2$ and corresponding involutive isometry $\sigma_2$, and a harmonic map $h : \underline{M}^n \to X$ such that (i) $h(B_1) \subset B_2$ and Iii) $h_*(w) \in T_{h(q)} B_2^{\perp}$ for $w \in T_q B_1^{\perp}$ for all $q \in B_1$. Then, we have $h(\sigma_1(p)) = \sigma_2(h(p))$ for all $p \in \underline{M}^n$.*

Proof: We will prove the theorem by an induction on the complex dimension n. For n =1, we can apply Theorem 4.1 directly with H trivial. Assuming that it is true for complex dimension n-1. Using the notation of Definition 6.1A, we set H = = $B_1 \cup \tilde{M}_{n-1}$ and apply Theorem 4.1 to complete the proof. □

In particular, let X be the real line **R**, $B_2$ be the point {0} and define $\sigma_2$ be $\sigma_2$ (y) = -y  for all  y ∈ R, then we have the following corolllary.

**Corollary 6.2.** *Keeping the notion of Theorem6.2, let  h  be a real harmonic function defined on  that vanishes on  then  $h(\sigma_1(p))$ = -h(p) for all p ∈ $\underline{M}^n$*

The following theorem is a direct consequence of Theorem 6.2, Lemma 6.1 and Remark 6.1 or Theorem 6.1 with a simple inductive argument.

**Theorem 6.3.** *Let f : M → N be a holomorphic map between the two Hermitian symmetric spaces, M is of complex dimension n,  $B_1$ be a recursive real form of M  and $B_2$ be a real form  N with involutive isometries $\sigma_1$  and  $\sigma_2$  (resp.) such that  $f(B_1) \subset B_2$ . Then, we have $f(\sigma_1(p))$ = $\sigma_2$ (f(p)) for all p ∈ M.*

**Corollary 6.3.** Let $\underline{M}^n$ an n dimensional n-dimensional Hermitian symmetric spaces and B is an recursive real form of  M with associated involutive anti-holomorphic isometry σ. Let f be a holomorphic function defined on M that takes on real values on the associated real form B, then f satisfies the following reflection principle

$$f(\sigma(p)) = \overline{f(p)} \quad \text{for all } p \in M^n .$$

Proof:

When $M^n$ = $\mathbf{C}^n$, a holomorphic function is considered as  a holomorphic map into the Gaussian plane. For the cases, it is considered as a holomorphic map into the Poincare unit disk.

We will use here the fact that if  a function  f  is holomorphic on $M^n$ and hence harmonic on $M^n$,  then  f is also holomorphic on  $\tilde{M}_k$  hence also harmonic on $\tilde{M}_k$ . □

For n=1, Corollary 6.2 and Corollary 6.3 specialize to the classical Schwarz reflection principle for harmonic and holomorphic functions defined on the complex plane. Even in the special case of a holomorphic function  f  defined on the complex plane **C** that takes on real values on the real line, the proof used in Theorem 6.2 is different from other commonly used proofs for such a reflection principle. In fact, the method of used herein allows us to prove the following more general fact.

**Theorem 6.3.** *Let  f  be a meromorphic function defined on  **C**  that takes on real  values, with possibly poles,  along the real line, then*

$$f(\bar{z}) = \overline{f(z)} ,$$

*for  z ∈ C for which f is defined.*

Proof: The function  f  defines an harmonic map of **C**  into itself. Define a new meromorphic function  g on **C** by

$$g(z) = \overline{f(\bar{z})}.$$

Let $p$ be a point on the real line of $\mathbb{C}$ at which $f$ is holomorphic and $U$ be a neighborhood of $p$ that contains no poles of $f$. Applying the proof of Theorem 6.2, we can show that

$$g(z) = f(z),$$

for all $z \in U$. Theorem 6.3 now follows from the uniqueness in the analytic continuation of a meromorphic function. □

**Remark 6.2.** The above proof can be applied to only to recursive real forms of non-compact Hermitian symmetric spaces. It is not clear that if the reflection principle for holomorphic functions is also true for other real forms of non-compact Hermitian symmetric spaces.

**Remark 6.3.** In this paper, we study the generalization of the Schwarz reflection principle mostly from a differential geometric perspective. Many deep and beautiful results for higher dimensional generalization of the Schwarz reflection principle from the perspectives of PDE and analysis can be found for examples in [BL] and [BR] for real hypersurfaces of complex spaces. It is not clear to the author that how such similar or related results could be proved for the recursive real forms.

# Appendix

## B. Real forms of Classical bounded symmetric domains.

| Type | dim$_C$ | Ambient space | Real forms |
|---|---|---|---|
| A III | pq | $\{Z \in M_{p,q}(C): {}^tZ\bar{Z} < I_q\} =$ SU(p, q)/S(U$_p$ x U$_q$) <br> p = q = 1 or p < q <br> p, q not both even <br> p, q both even <br><br> p = q > 1 <br> p odd <br> p even | <br><br><br> SO(p, q)/SO(p) x SO(q) <br> SO(p, q)/SO(p) x SO(q); <br> Sp(p/2, q/2)/Sp(p/2) x Sp(q/2) <br><br> SO(p, p)/SO(p) x SO(p); SL(p, C) x R <br> SO(p, p)/SO(p) x SO(p); SL(p, C) x R; <br> Sp(p/2, p/2)/Sp(p/2) x Sp(p/2) |
| D III | n(n-1)/2 | $\{Z \in M_{n,n}(C): {}^tZ = -Z, {}^tZ\bar{Z} < I_n\} =$ SO*(2n)/U(n) <br> n odd <br> n even | <br><br> SO(n, C) <br> SO(n, C); [SU*(n)/Sp(n/2)] x R |
| BD I (q=2) | p | $\{X \in M_{2,p}(R): X{}^tX < I_2\} =$ SO(2, p)/SO(p) x SO(2) | [SO(1, k)/SO(1) x SO(k)] x [SO(1, p – k + 1)/SO(1) x SO(p-k+1)], $0 \leq k \leq [p/2]$ |
| C I | n(n+1)/2 | $\{Z \in M_{n,n}(C): {}^tZ = -Z, {}^tZ\bar{Z} < I_n\} =$ Sp(n, R)/U(n) <br> n odd <br> n even | <br><br> [SL(n, R)/SO(n)] x R <br> [SL(n, R)/SO(n)] x R; Sp(n/2, C) |

E-mail Addresses:  domleung_cal70@berkeley.edu